\documentclass[8pt,a4,twoside,rm]{article}
\usepackage{amssymb}
 \usepackage{amsthm}
\usepackage{color}
\usepackage{pstricks}
\usepackage{hyperref}
\usepackage{amsmath}
\usepackage[left=2cm,top=4cm,right=2cm,bottom=4cm]{geometry}
\usepackage{amsmath, amsthm, amssymb}
\usepackage{graphicx}
\usepackage{amssymb}
\usepackage{fancyhdr}
\usepackage{titlesec}
\usepackage{youngtab}
\usepackage{color}
\usepackage{hyperref}
\usepackage{pstricks}
\usepackage{pst-all}

\newtheorem{thmx}{Theorem}

\newtheorem{teo}{Theorem}[section]
\newtheorem{lem}[teo]{Lemma}

\newtheorem{pro}[teo]{Proposition}
\newtheorem{defi}[teo]{Definition}
\newtheorem{cor}[teo]{Corollary}

\newtheorem{rem}[teo]{Remark}

\newtheorem{cho}[teo]{Choices}

\newenvironment{dem}{\noindent \textit{Proof:} }{\quad \hfill $\square$}

\newcommand{\hm}{\mbox{Hom}}
\newcommand{\sub}[1]{ \underline{#1} }

\newcommand{\K}{\mathbb R}
\newcommand{\R}{\mathbb R}
\newcommand{\N}{\mathbb N}
\newcommand{\Z}{\mathbb Z}
\newcommand{\bs}[1]{ BS(\sub{#1})}
\newcommand{\bsa}[1]{ BS(#1)}
\newcommand{\al}[1]{ A_{\underline{#1}}}
\newcommand{\cell}[2]{ \Delta_{\underline{#1}}(#2) }
\newcommand{\simple}[2]{ D_{\underline{#1}}(#2) }
\newcommand{\simpleset}[1]{ \Lambda_{0}(\underline{#1}) }
\newcommand{\gdn}[3]{ d_{\sub{#1}}(#2,#3) }
\newcommand{\res}[1]{ \mbox{Res}(#1) }
\newcommand{\gd}{\dim_{v}}
\newcommand{\overlim}[1]{{\buildrel{#1}\over\longrightarrow\;}}
\newcommand{\fre}[1]{ \sub{\mathbf{#1}}}

\def\levelonetree{
\pspicture(-2,-2)(8,4)
\psline{->}(3,3)(-1,-2)
\psline{->}(3,3)(7,-2)
\rput{0}(3,3.4){$\bsa{s_1\ldots s_n}$}
\rput{0}(-1,-2.4){$\bsa{s_1\ldots s_{n-1}}$}
\rput{0}(7,-2.4){$\bsa{s_1\ldots s_n}$}
\rput{0}(-1,1){$ \mathbb{I}^{n-1} \otimes m_{s_n}$}
\rput{0}(6,1){$\mathbb{I}$}
\endpspicture
}

\def\leveltwotree{
\pspicture(-3.5,-2)(9.5,4)
\rput(3,3.4){$\bsa{s_{1}\ldots s_{n-k+1}} \bsa{\fre{x}} $}
\rput(-3,-2.4){$\bsa{s_{1}\ldots s_{n-k}} \bsa{\fre{x}} $}
\rput(9,-2.4){$\bsa{s_{1}\ldots s_{n-k}} \bsa{\fre{s_{n-k+1}x}} $}
\rput(-3,1){$ \mathbb{I}^{n-k} \otimes m_{s_{n-k+1}} \otimes \mathbb{I}_{\fre{x}} $}
\rput(10,1){$\mathbb{I}^{n-k} \otimes F(\sub{s_{n-k+1} x},\fre{s_{n-k+1}x} )$}
\psline{->}(3,3)(-3,-2)
\psline{->}(3,3)(9,-2)
\endpspicture
}

\def\levelthreetree{
\pspicture(-6,-3)(6,14.5)
\rput(0,14.5){$ \bsa{s_{1}\ldots s_{n-k+1}} \bsa{\fre{x}}$}
\rput(0,11){$ \bsa{s_{1}\ldots s_{n-k+1}} \bsa{\sub{x}_{s_{n-k+1}}}$}
\rput(0,7){$\bsa{s_1\ldots s_{n-k}} \bsa{\sub{x}_{s_{n-k+1}}}  $ }
\rput(-6,2){$\bsa{s_{1}\ldots s_{n-k}}  \bsa{\sub{s_{n-k+1}x}}  $}
\rput(6,2){$\bsa{s_{1}\ldots s_{n-k}}  \bsa{\sub{x}_{s_{n-k+1}}}$}
\rput(-6,-2){$\bsa{s_{1}\ldots s_{n-k}}  \bsa{\fre{s_{n-k+1}x}}  $}
\rput(6,-2){$\bsa{s_{1}\ldots s_{n-k}}  \bsa{\fre{x}}$}
\psline{->}(0,14)(0,11.5)
\psline{->}(0,10.5)(0,7.5)
\psline{->}(0,6.5)(-4,2.5)
\psline{->}(0,6.5)(4,2.5)
\psline{->}(-6,1.5)(-6,-1.5)
\psline{->}(6,1.5)(6,-1.5)
\rput(3.3,12.5){$ \mathbb{I}^{n-k+1} \otimes F(\fre{x},\sub{x}_{s_{n-k+1}}) $}
\rput(3.2,9){$ \mathbb{I}^{n-k}\otimes j_{s_{n-k+1}} \otimes \mathbb{I}^{l(x)-1}$}
\rput(-5,5){$ \mathbb{I}^{n-k}\otimes m_{s_{n-k+1}} \otimes \mathbb{I}^{l(x)-1}$}
\rput(3,5){$ \mathbb{I}$}
\rput(-10.5,0){$ \mathbb{I}^{n-k} \otimes F(\sub{s_{n-k+1}x},\fre{s_{n-k+1}x}) $}
\rput(9,0){$ \mathbb{I}^{n-k} \otimes F(\sub{x}_{s_{n-k+1}},\fre{x}) $}
\endpspicture }

\def\newleavesfromoldleaves{
\pspicture(-6,0)(9,8)
\rput(-6,2){$\bsa{\fre{x}}$}
\rput(-2,2){$\bsa{\fre{x'}}$}
\rput(-4,5){$ B_s \otimes_{R} \bsa{\fre{x}}$}
\rput(-4,8){$\bs{w} =B_{s} \otimes_{R} \bs{w'}$}
\psline{->}(-4,7.7)(-4,5.3)
\psline{->}(-4,4.7)(-6,2.3)
\psline{->}(-4,4.7)(-2,2.3)
\rput(-3,6.5){$\mathbb{I}_{s} \otimes l$}
\rput(-7,3.5){$j_{s} \otimes \mathbb{I}_{\fre{x'}}$}
\rput(-1,3.5){$m_sj_{s} \otimes \mathbb{I}_{\fre{x'}}$}
\rput(5,2){$\bsa{\fre{x}}$}
\rput(9,2){$\bsa{\fre{x'}}$}
\rput(7,5){$ B_s \otimes_{R} \bsa{\fre{x}}$}
\rput(7,8){$\bs{w} =B_{s} \otimes_{R} \bs{w'}$}
\psline{->}(7,7.7)(7,5.3)
\psline{->}(7,4.7)(9,2.3)
\psline{->}(7,4.7)(5,2.3)
\rput(8,6.5){$\mathbb{I}_{s} \otimes l$}
\rput(4,3.5){$m_s \otimes \mathbb{I}_{\fre{x}}$}
\rput(10,3.5){$\mathbb{I}_s \otimes \mathbb{I}_{\fre{x}}$}
\rput(-4,1){$x'<x$}
\rput(7,1){$x'>x$}
\endpspicture
	}

\numberwithin{equation}{section}

\titleformat{\section}[hang]{\sc}{\thesection.}{0.5cm}{\filcenter}
\titleformat{\subsection}[hang]{\bf}{\thesubsection.}{0.2cm}{\filright}

\begin{document}

\title{ \normalsize {\bf \sc Categorification of a recursive formula for Kazhdan-Lusztig polynomials}}
\author{\sc  David Plaza }
\date{ }
\maketitle

\begin{center}
Universidad de Chile,\\
Facultad de Ciencias,\\
Santiago, Chile
\end{center}

\begin{center}
Email: davidricardoplaza@gmail.com
\end{center}

\begin{abstract}
We obtain explicit branching rules for graded cell modules and graded simple modules over the endomorphism algebra of a Bott-Samelson bimodule. These rules allow us to categorify a well-known recursive formula for Kazhdan-Lusztig polynomials.
\end{abstract}







\section{Introduction}   \label{section intro}
Let $(W,S)$ be a Coxeter system.  Let $l:W \rightarrow \Z_{\geq 0}$ be the corresponding length function. We denote by $\leq $ to the usual Bruhat order on $W$.   Kazhdan and Lusztig (KL) introduced a family of polynomials with integer coefficients that are indexed using pairs of elements of $W$ \cite{kazhdan1979representations}. These polynomials are now known as KL polynomials and can be defined in multiple equivalent ways. The definition provided here best suits our purposes and involves a recursive process. For  $x,w \in W$, we define the KL polynomial, denoted $P_{x,w} \in \Z[q]$, as follows. First, we set
\begin{equation}   \label{eq intro recursive}
P_{x,w}=0 \mbox{ if } x\not \leq w, \mbox{ and } P_{x,w}=1 \mbox{ if } x=w.
\end{equation}

Let $e$ be the identity element of W. Recall that $e \leq x$, for all $x \in W$. Therefore, $P_{x,e}$ is determined by  (\ref{eq intro recursive}), for all $x\in W$. Now, fix $w \in W$ and assume that $P_{u,v}$ has been defined, for all $u\in W$ and $v < w$.  For  $x<w$ and  $s\in S$, satisfying $l(sw)<l(w)$, we define
\begin{equation}
c=c_{s}(x):=\left\{
    \begin{array}{ll}
      0, & \mbox{ if }  x <sx ; \\
      1, & \mbox{ if }  x>sx.
    \end{array}
  \right.
\end{equation}
Then, $P_{x,w}$ is defined as
\begin{equation}  \label{equation introduction}
P_{x,w}= q^{c-1}P_{sx,sw} +q^{c}P_{x,sw}- \sum_{sz< z < sw} \mu (z,sw) q^{(l(z)-l(w))/2}P_{x,z},
\end{equation}
where $\mu (z,sw)$ is the coefficient of $q^{(l(sw)-l(z)-1)/2}$ in $P_{z,sw}$. It was proven in \cite{humphreys1992reflection} that this definition is equivalent to the original definition given by Kazhdan and Lusztig in \cite{kazhdan1979representations}.

\medskip
The main result in this paper is a categorification of (\ref{equation introduction}). The rest of this section provides a precise explication about  what is meant by a categorification of (\ref{equation introduction}). Roughly speaking, we refer to a categorification of (\ref{equation introduction}) as the process that gives  (\ref{equation introduction}) using the category of Soergel bimodules.

\medskip
For every Coxeter system $(W,S)$, Soergel \cite{soergel1992combinatorics} constructed a category, $\mathcal{B}=\mathcal{B}(W,S,V,\K)$, of $\Z$-graded bimodules over a polynomial ring $R$ with coefficients in the real numbers $\R$. It depends on $(W,S)$ and on a finite-dimensional $\R$-representation $V$ of $W$. Soergel proved that (up to a degree shift) $W$ parameterizes the set of indecomposable objects in $\mathcal{B}$. We let $B_{w}\in \mathcal{B}$ denote the corresponding indecomposable object. Soergel \cite{soergel2007kazhdan} also proved that $\mathcal{B}$ is a categorification of the Hecke algebra $\mathcal{H}$
of $W$.  This means that there exists an algebra isomorphism $\eta: [\mathcal{B}] \rightarrow \mathcal{H}$  between the split Grothendieck group $[\mathcal{B}]$ of $\mathcal{B}$ and $\mathcal{H}$. Soergel proposed the following conjecture, which came to be known as Soergel's conjecture:
\begin{equation}
\eta([B_w]) = \sub{H}_{w},
\end{equation}
where $\{\sub{H}_w\}_{w\in W}$ is the KL basis of $\mathcal{H}$. In 2014, Elias and Williamson \cite{elias2012hodge} proved this conjecture.
We should mention that, by applying Soergel's previous work, their results provide a proof for the longstanding KL \emph{positivity conjecture} for every Coxeter system. In a previous publication  \cite{plaza2014graded}, the author proved that Soergel's conjecture also implies a related conjecture for KL polynomials, known as the \emph{monotonicity conjecture}.

\medskip
The category of Soergel bimodules can be obtained as the graded Karoubi envelope of another category, $\mathbb{BS}\mbox{Bim}$, of $\Z$-graded $(R,R)$-bimodules, known as Bott-Samelson bimodules. For each \emph{expression} $\sub{w}= (s_1,s_2, \ldots , s_k) \in S^k$, we have a Bott-Samelson bimodule $\bs{w}$ (for details see Section \ref{section : the category of soergel bimodules}). Let $R^{+}$ be the ideal of $R$ generated by homogeneous elements of nonzero degree. As
$\R \cong R/ R^{+}$ we can view  $\R$ as a left $R$-module. The endomorphism ring $\mbox{End}(\bs{w})$ of a Bott-Samelson bimodule $\bs{w}$ can be naturally equipped with a $(R,R)$-bimodule structure. Therefore, we can define the $\R$-algebra $\al{w}:= \mbox{End}(\bs{w}) \otimes_{R} \R$. The algebra $\al{w}$ is a graded cellular $\R$-algebra in the sense of Hu and Mathas \cite{hu2010graded}, with a graded cellular basis, the \emph{double leaves basis}. This basis was defined by Libedinsky in \cite{libedinsky2015light} and it generalizes his previous construction of the \emph{light leaves basis} \cite{libedinsky2008categorie}.

\medskip
The existence of a graded cellular basis for $\al{w}$ allows us to define graded (right) cell modules and graded (right) simple modules, as well as graded decomposition numbers. If $\sub{w}$ is a reduced expression of an element $w \in W$, then the set of graded cell modules for $\al{w}$ is parameterized by

\begin{equation}
\Lambda(\sub{w}):=\{ x \in W \mbox{ }| \mbox{ }x \leq w \}.
\end{equation}

Given $x \in  \Lambda(\sub{w})$, we use $\cell{w}{x}$ to denote the corresponding graded cell module. Graded simple modules are parameterized by a subset $\Lambda_{0}(\sub{w})$ of $\Lambda(\sub{w})$. For $y \in \Lambda_{0}(\sub{w}) $, we denote by $\simple{w}{y}$  the corresponding graded simple module. The set $\{  \simple{w}{y} \mbox{ } | \mbox{ } y \in \simpleset{w} \}$  is a complete set (up to a degree shift) of pairwise non-isomorphic graded simple right $\al{w}$-modules. We also let $\gdn{w}{x}{y} \in \Z[v,v^{-1}]$ denote the corresponding graded decomposition number. We recall that $\gdn{w}{x}{y} $ counts the number of times $\simple{w}{y}$ appears in a graded composition series of $\cell{w}{x}$. Because Soergel's conjecture proved to be valid, these graded decomposition numbers coincide with the KL polynomials when the former are defined and the latter are suitably normalized (see Proposition \ref{equality graded decompostion numbers with KL}).

\medskip
For the rest of this section we fix an element $w\in W$ and a reduced expression $\sub{w}=(s_1,s_2,\ldots, s_{k})$ of $w$. The expression $\sub{w'}=(s_2, \ldots ,s_{k})$ is a reduced expression for $w':=s_1w \in W$. The algebra $\al{w'}$ can be naturally embedded into $\al{w}$. Accordingly, each $\al{w}$-module can be considered as a $\al{w'}$-module by restriction.  Given a graded $\al{w}$-module $M$, we let $\mbox{Res} (M)$ represent its restriction. In particular, we can restrict  graded cell modules and graded simple modules of $\al{w}$. In this paper, we obtain explicit graded branching rules for these modules. In other words, we show how the graded cell modules and the graded simple modules of $\al{w}$ decompose (when they are viewed as $\al{w'}$-modules) in terms of graded cell modules and graded simple modules of $\al{w'}$, respectively.

\medskip
A concise way to express these branching rules is in terms of graded Grothendieck groups. Let $A$ be a finite dimensional $\Z$-graded algebra. The \emph{graded Grothendieck group} $\mathcal{G}(A)$ of $A$ is the  $\Z[v,v^{-1}]$-module generated by symbols $[M]$, where $M$ runs over the finite dimensional graded right $A$-modules, with relations:
\begin{enumerate}
\item   $[M\langle k \rangle]= v^{k} [M]$, for all graded right $A$-module $M$ and  $k\in \Z$;
\item   $[M]=[M']+[M'']$ if there exists a  short exact sequence $   0\rightarrow M' \rightarrow M \rightarrow M'' \rightarrow 0  $ of  graded right $A$-modules.
\end{enumerate}

In particular, for $A= \al{w}$ (resp. $A=\al{w'}$), we write $\mathcal{G}_{\sub{w}}:=\mathcal{G}(\al{w})$ (resp.
$\mathcal{G}_{\sub{w'}}:=\mathcal{G}(\al{w'}))$. Using this terminology, the graded branching rules can be expressed as follows:

\begin{thmx}   \label{teoprincipalintro}
 There exists a homomorphism of $\Z[v,v^{-1}]$-modules $\mbox{Res}: \mathcal{G}_{\sub{w}} \rightarrow \mathcal{G}_{\sub{w'}}$ determined by
\begin{equation}
 \res{[M]} = [\res{M}],
\end{equation}
for every right graded $\al{w}$-module $M$. The image of the class of the cell module $\cell{w}{x}$ under this homomorphism is
\begin{equation}
     \mbox{Res}([\cell{w}{x}] )=  \left\{
                                  \begin{array}{rl}
                                    v^{-1} [\cell{w'}{x}] + [\cell{w'}{s_1x}], & \mbox{if } s_1x<x ; \\
                                     v [\cell{w'}{x}] + [\cell{w'}{s_1x}], & \mbox{if }  s_1x>x.
                                  \end{array}
                                \right.
\end{equation}
On the other hand, if $y\in \simpleset{w}$, then the restriction of the simple module $\simple{w}{y}$ is given by
\begin{equation}
\mbox{Res}([\simple{w}{y}]) = \sum_{u\in \simpleset{w'}} h_{s_1,u}^{y} [\simple{w'}{u}],
\end{equation}
where $ h_{s_1,u}^{y}\in \Z [v,v^{-1}]$ denotes the coefficient of $\sub{H}_y$ in the expansion of $\sub{H}_{s_1}\sub{H}_{u}$ in terms of the KL basis of $\mathcal{H}$.
\end{thmx}

From the definition of Grothendieck group, we know that  $\mathcal{C}_{\sub{w}}:= \{ [\simple{w}{y}] \mbox{ } | \mbox{ } y\in \simpleset{w}  \}$  is a
$\Z [ v,v^{-1} ]$-basis of $\mathcal{G}_{\sub{w}}$. Furthermore, by definition of graded decomposition numbers, for $x\leq w$  the expansion of the class $[  \cell{w}{x} ]$ in terms of the basis $\mathcal{C}_{\sub{w}}$ is given by
\begin{equation}   \label{equation intro cate}
[\cell{w}{x}]=  \sum_{y\in \simpleset{w}} \gdn{w}{x}{y} [\simple{w}{y}].
\end{equation}

We can now apply the homomorphism $\mbox{Res}$ to the equation (\ref{equation intro cate}). Theorem \ref{teoprincipalintro} provides two different ways to calculate this value. This allows us to obtain the following equation in $\mathcal{G}_{\sub{w'}}$.

\begin{equation}   \label{equation  intro cate dos}
\sum_{u\in \simpleset{w'}} \left(\sum_{y\in \simpleset{w}}\gdn{w}{x}{y} h_{s_1,u}^{y} \right)[\simple{w'}{u}]= \left\{
                                  \begin{array}{rl}
                                    v^{-1} [\cell{w'}{x}] + [\cell{w'}{s_1x}], & \mbox{if } s_1x<x ; \\
                                     v [\cell{w'}{x}] + [\cell{w'}{s_1x}], & \mbox{if }  s_1x>x.
                                  \end{array}
                                \right.
\end{equation}

Equating coefficients of the basis elements $[\simple{w'}{u}]\in \mathcal{C}_{\sub{w'}}  $ yields

\begin{equation}   \label{equation  intro cate tres}
\sum_{y\in \simpleset{w}}\gdn{w}{x}{y} h_{s_1,u}^{y} = \left\{
                                  \begin{array}{rl}
                                    v^{-1}  \gdn{w'}{x}{u} + \gdn{w'}{s_1x}{u}, & \mbox{if } s_1x<x ; \\
                                     v \gdn{w'}{x}{u} + \gdn{w'}{s_1x}{u}, & \mbox{if }  s_1x>x,
                                  \end{array}
                                \right.
\end{equation}
for all $u\in \simpleset{w'}$. As we already mentioned, the graded decomposition numbers coincide with (a normalized version of) the KL-polynomials. Accordingly, equation (\ref{equation  intro cate tres}) provides a relation between  KL-polynomials.  In general, determining which elements belong to $\simpleset{w'}$ is a hard task. However, we can straightforwardly note that $w' \in \simpleset{w'}$. Equation (\ref{equation  intro cate tres}) corresponding to the element $w'$ is a normalized version of (\ref{equation introduction}) (as we will see in Section 5). Summing up, we have produced  (\ref{equation introduction}) by studying the category of Soergel bimodules. We consider this procedure to be a categorification of (\ref{equation introduction}).

\medskip
The rest of this paper is organized as follows. In the next section, given an arbitrary Coxeter system, we introduce its corresponding Hecke algebra and Soergel bimodule category. We also show how to obtain a normalized version of (\ref{equation introduction}). Section 3 reviews Libedinsky's construction of the double leaves basis. We recall that the double leaves basis is a graded cellular basis in Section 4 and then prove that the graded decomposition numbers coincide with the KL polynomials. Finally, in Section 5, we obtain graded branching rules for the graded cellular and graded simple modules. We conclude this section with a categorification of a normalized version of (\ref{equation introduction}) using these branching rules.

\section{Preliminaries}
In this section, given an arbitrary Coxeter group, we introduce its Hecke algebra and its corresponding category of Soergel bimodules. We discuss the relation of these subjects with the KL polynomials and also present a process for obtaining a normalized version of (\ref{equation introduction}).

\subsection{Hecke algebras and Kazhdan-Lusztig polynomials} \label{hecke and kl poly}

Let $(W,S)$ be a Coxeter system with length function $l: W \rightarrow \N$, and let $e\in W$ denote the identity. We denote the order of $st \in W$ by $m_{st}=\{1,2, \ldots , \infty \}$, for all $s,t\in S$. The Hecke algebra $\mathcal{H}=\mathcal{H}(W,S)$ of $(W,S)$ is the $\Z[v,v^{-1}]$-algebra with generators $\{ H_{s} \mbox{ } | \mbox{ } s\in S \}$ and relations
\begin{equation} \label{relations hecke algebra one}
H_{s}^{2}= (v^{-1}-v)H_s+1    \mbox{ for all } s\in S,
\end{equation}

\begin{equation} \label{relations hecke algebra two}
\underbrace{H_sH_tH_s \ldots}_{m_{st}\mbox{ times} }  =  \underbrace{H_tH_sH_t \ldots }_{m_{st} \mbox{ times}  }    \mbox{ for all } s\neq t \in S.
\end{equation}

If $m_{st}=\infty$, then the \emph{braid relation} (\ref{relations hecke algebra two}) is omitted. Throughout this paper, we use an underlined Roman letter to denote a finite sequence of elements in $S$, and we call such a  sequence an \emph{expression}. If $\underline{w} = s_1s_2 \ldots s_k$ is an expression, then its length is $l(\underline{w})=k$. Omission of the underlining in an expression denotes the respective product in $W$. A expression $\underline{w} = s_1s_2 \ldots s_k$ is  \emph{reduced} if $l(\underline{w})= l(w)$, where $l(w)$ is the length function of $W$. Differentiating between expressions and elements in $W$ is important because many of the concepts defined in this paper depend on expressions rather than elements in $W$.

\medskip
Given an element $w\in W$ and a reduced expression $\underline{w}=s_1s_2 \ldots s_k$ of $w$, we define $H_{w}=H_{s_{1}} H_{s_2} \ldots H_{s_k} \in \mathcal{H}$. It follows from (\ref{relations hecke algebra two}) that $H_{w}$ is well defined, i.e., $H_w$ does not depend on the choice of a reduced expression for $w$. The set $\{H_{w} \mbox{ }  | \mbox{ } w\in W \}$ is a basis of $\mathcal{H}$ as a $\Z[v,v^{-1}]$-module. We refer to this basis as the \emph{standard basis} of $\mathcal{H}$. It follows from (\ref{relations hecke algebra one}) that $H_{s}$ is invertible in $\mathcal{H}$, for all $s\in S$. Thus, $H_{w}$ is also
invertible for every $w\in W$. Moreover, there exists a unique ring involution $d: \mathcal{H} \rightarrow \mathcal{H}$ determined by $d(v)=v^{-1}$ and $d(H_{w})=H_{w^{-1}}^{-1}$, for all $w\in W$.

\begin{teo} (Kazhdan-Lusztig \cite{kazhdan1979representations})  \label{teo kazhdan-lusztig}
There exists a unique basis $\{   \underline{H}_{w} \}_{ w\in W}$ of $\mathcal{H}$ as a $\Z[v,v^{-1}]$-module that satisfies
\begin{equation} \label{condition kl}
 d({\sub{H}}_{w}) = \sub{H}_{w}, \qquad \mbox{ and } \qquad
\underline{H}_{w} = H_{w} + \sum_{x < w} h_{x,w} H_{x},
\end{equation}
in which $h_{x,w} \in v\Z[v]$, and $<$ denotes the usual Bruhat order on $W$.
\end{teo}

The set $\{   \underline{H}_{w} \mbox{ }| \mbox{ } w\in W\}$ is called the \emph{Kazhdan-Lusztig basis} of $\mathcal{H}$, and the polynomials $h_{x,w} \in \Z[v]$ are known as the \emph{Kazhdan-Lusztig polynomials}.

\begin{rem} \rm \label{variables}
We should note that this follows the normalization provided by Soergel in \cite{soergel1997kazhdan} rather than the original normalization provided by Kazhdan and Lusztig in \cite{kazhdan1979representations}, which was used in the introduction. The two normalizations are related by  $q = v^{-2}$, and the original Kazhdan-Lusztig polynomials $P_{x, w}(q)\in \mathbb{Z}[q]$ can be recovered from the Kazhdan-Lusztig polynomials $h_{x, w}(v)\in \mathbb{Z}[v]$ using
\begin{equation}
h_{x,w}(v) = v^{l(w) - l(x)}P_{x, w}(v^{-2}).
\end{equation}
\end{rem}

Theorem \ref{teo kazhdan-lusztig} provides an algorithm that can be used to calculate (inductively on the Bruhat order) the KL basis and the KL polynomials. The first KL basis elements can be easily defined as
\begin{equation}  \label{first kl basis}
\sub{H}_{e} := H_e  \qquad \mbox{ and } \qquad  \sub{H}_{s} =: H_{s} +v, \qquad  \mbox{ for all } s\in S.
\end{equation}
In order to continue this process, we must know the action of $\sub{H}_s$ on the standard basis of $\mathcal{H}$. This action is given by
\begin{equation}  \label{equation action Hs}
\sub{H}_sH_x = \left\{
                 \begin{array}{ll}
                   H_{sx}+vH_{x}, & \mbox{ if } sx>x; \\
                  H_{sx}+v^{-1}H_{x}, & \mbox{ if } sx<x,
                 \end{array}
               \right.
\end{equation}
for all $x\in W$. Let us return to the calculation of the elements of the KL basis. Fix $u \in W$ and suppose $\sub{H}_{z}$ has been calculated, for all $z\leq u$. We choose $s\in S$ with $su > u$ and write $\sub{H}_s \sub{H}_{u}$ in terms of the standard basis of $\mathcal{H}$:
\begin{equation}  \label{alg kl poly}
\sub{H}_s \sub{H}_{u} = H_{su} +\sum_{z < su} g_{z} H_{z},
\end{equation}
for some $g_{z} \in \Z[v,v^{-1}]$. By applying (\ref{equation action Hs}) and from the fact that the KL polynomials have coefficients in $v\Z[v]$, we conclude that $g_{z} \in \Z[v]$. Hence, it is easy to see that
\begin{equation} \label{alg kl poly dos}
\sub{H}_{su} = \sub{H}_{s}\sub{H}_{u} - \sum_{z < su} g_{z}(0) \sub{H}_{z}
\end{equation}
because the right side of (\ref{alg kl poly dos}) satisfies the conditions (\ref{condition kl}) provided in Theorem \ref{teo kazhdan-lusztig}. Summing up, $\sub{H}_{su}$ can be calculated using our previous knowledge of all the elements $\sub{H}_{z}$ of the KL basis in which $z\leq u$. We should note that this algorithm is far from efficient because a substantial amount of information is necessary for each new step in the algorithm.

\medskip
On the other hand, equation (\ref{alg kl poly dos}) also provides an inductive procedure that can be used to compute the KL polynomials. Concretely, let $w = su$ and choose $x\in W$ with $x < w$. Equating the coefficients of $H_x$ on both sides of (\ref{alg kl poly dos}) yields
\begin{equation}  \label{alg kl poly cuatro}
h_{x,w} = \left\{
            \begin{array}{rl}
        \displaystyle    v h_{x,sw} + h_{sx,sw} - \sum_{z < w} g_{z}(0) h_{x,z}  , & \mbox{ if } x<sx ; \\
        \displaystyle     v^{-1} h_{x,sw} + h_{sx,sw}  - \sum_{z < w} g_{z}(0) h_{x,z}, & \mbox{ if } x>sx.
            \end{array}
          \right.
\end{equation}

In addition, (\ref{equation action Hs}) implies that $g_{z}(0)\neq 0$ only if $sz < z < sw$. Furthermore, equation (\ref{alg kl poly}), which defines the polynomials $g_{z}$, allows us to conclude that $g_{z}(0)$ equals the linear coefficient of the KL polynomial $h_{z,sw}$. We use $\mu (z,sw)$ to denote this value. Thus, we can rewrite (\ref{alg kl poly cuatro}) as

\begin{equation}  \label{alg kl poly cinco}
h_{x,w} = \left\{
            \begin{array}{rl}
        \displaystyle    v h_{x,sw} + h_{sx,sw} - \sum_{sz< z < sw}   \mu (z,sw)  h_{x,z}  , & \mbox{ if } x<sx ; \\
        \displaystyle     v^{-1} h_{x,sw} + h_{sx,sw}  - \sum_{sz <z < sw} \mu (z,sw)  h_{x,z}, & \mbox{ if } x>sx.
            \end{array}
          \right.
\end{equation}

Note that we have already used the symbol $\mu (z,sw)$ in (\ref{equation introduction}) to denote the coefficient of $q^{(l(sw)-l(z)-1)/2}$ in the \emph{original} KL polynomial $P_{z,sw}$. However, Remark \ref{variables} makes it clear that the linear coefficient of $h_{z,sw}$ coincides with the coefficient of $q^{(l(sw)-l(z)-1)/2}$ in $P_{z,sw}$. This remark further indicates that (\ref{alg kl poly cinco}) is a normalized version of (\ref{equation introduction}). As mentioned in the introduction, the main result in this paper is a categorification of (\ref{equation introduction}), or, equivalently, a categorification of (\ref{alg kl poly cinco}). In actuality, we categorify (\ref{alg kl poly cinco}) rather than (\ref{equation introduction}) because the former is best suited for our purposes.

\subsection{The category of Soergel bimodules}  \label{section : the category of soergel bimodules}

Recall that $(W,S)$ is an arbitrary Coxeter group. In order to introduce the category of Soergel bimodules, we must first fix a representation $V$ of $W$, which should satisfy certain technical requirements.

\begin{defi}  \rm
A \emph{reflection faithful} representation $V$ of $W$ over $\K$ is a finite-dimensional $\K$-representation of $W$ that satisfies the following conditions.
\begin{enumerate}
\item  The representation $V$ is faithful.
\item  If $V^{w}$ is the set of elements in $V$ fixed by $w \in W$, then $V^{w}$ has codimension $1$ if and only if $w$ is conjugated to a simple reflection.
\end{enumerate}
\end{defi}

In \cite{soergel2007kazhdan}, Soergel demonstrated that such a representation exists for arbitrary Coxeter groups. Therefore, we can fix (for the rest of the paper) a reflection-faithful representation $V$ of $W$ defined over $\R$. Let $R$ be the $\R$-algebra of regular functions on $V$ with the following grading:
\begin{equation}
R= \bigoplus_{i\in \Z_{\geq 0}} R_{i}\mbox{, with } R_2=V^{*} \mbox{ and }  R_{i}=0 \mbox{ if } i \mbox{ is odd.}
\end{equation}
We can think of the elements in $R$ as polynomials over $V^*$. We also define $R^{+}$ to be the subring of $R$ generated by all homogeneous elements having a nonzero degree. It is easy to see that $\R \cong R/ R^{+}$. This isomorphism will be frequently used in order to consider $\R$ as a left $R$-module.

\medskip
Given a graded $(R,R)$-bimodule $B$ and $k\in \Z$, we let $B(k)$ signify the graded $(R,R)$-bimodule that is obtained from $B$ by shifting the grading by $k$. In other words, if $B=\bigoplus_{i\in \Z} B_i$ is a graded bimodule, then $B(k)_{i} = B_{k+i}$ for all $i\in \Z$. There is a natural action of $W$ on $R$, which is induced by the action of $W$ on $V$. For $s\in S$, let $R^{s}$ be the subring of $R$ fixed by $s$. The graded $(R,R)$-bimodule $B_s$ is defined to be
\begin{equation}
B_s= R \otimes_{R^{s}} R  (1),
\end{equation}
where the left (resp. right) action of $R$ on $B_{s}$ is given by left (resp. right) multiplication. Given any expression $\underline{w}= s_{1}s_{2} \ldots s_{k}$, we define the \emph{Bott-Samelson bimodule}  $BS(\underline{w})$ to be
\begin{equation}
BS(\underline{w})= B_{s_1} \otimes_{R}  B_{s_{2}} \otimes_{R} \ldots \otimes_{R} B_{s_{k}}.
\end{equation}
We have the following $(R,R)$-bimodule isomorphism.
\begin{equation}
BS(\underline{w})\cong R \otimes_{R^{s_{1}}}  R  \otimes_{R^{s_{2}}}  \ldots \otimes_{R^{s_{k}}} R (k)
\end{equation}
 Therefore, we can write any element of this module as a sum of terms given by $k + 1$ polynomials in $R$, one in each slot, separated by the tensors. Denote by  $\mathbb{BS}\mbox{Bim}$  the category whose objects include all direct sums and grading shifts of  Bott-Samelson bimodules and whose morphisms are degree-preserving bimodule homomorphisms. For two objects $B, B' \in \mathbb{BS}\mbox{Bim}$, let $\mbox{Hom} (B,B')$ be the respective morphisms space. We also define
\begin{equation}
\mbox{Hom}^{\Z} (B,B')= \bigoplus_{k\in  \Z} \mbox{Hom} (B(k), B').
\end{equation}
Finally, we define the category of Soergel bimodules, $\mathbb{S}\mbox{Bim}$, as the Karoubi envelope of $\mathbb{BS}\mbox{Bim}$. In other words, the objects in $\mathbb{S}\mbox{Bim}$ are direct sums and graded shifts of direct summands of the objects in $\mathbb{BS}\mbox{Bim}$. In \cite[Theorem 6.16]{soergel2007kazhdan}, Soergel demonstrated that the indecomposable objects in $\mathbb{S}\mbox{Bim}$ are indexed using elements of $W\times \Z$. They are written as $B_{w}(k)$ for $w\in W$ and $k\in \Z$. An indecomposable bimodule $B_{w}$ is completely determined by one property: it appears as a direct summand of $\mathbb{BS}(\underline{w})$ for every reduced expression $\underline{w}$ of $w$, and it does not appear as a direct summand of $\mathbb{BS}(\underline{x})$ for any expression $\underline{x}$ with length less than $l(w)$.

\medskip
Denote by $[\mathbb{S}\mbox{Bim}]$  the split Grothendieck group of $\mathbb{S}\mbox{Bim}$. That is, $[\mathbb{S}\mbox{Bim}]$ is the abelian group generated by $[B]$ for every object $B \in \mathbb{S} \mbox{Bim}$, subject to the relation $[B]=[B_1]+[B_2]$ whenever $B \cong B_1\oplus B_{2}$. If we define $[B_1][B_{2}] = [B_1 \otimes B_2]$, then $[\mathbb{S}\mbox{Bim}]$ is equipped with a ring structure. Furthermore, $[\mathbb{S}\mbox{Bim}]$ has the structure of  $\Z[v,v^{-1}]$-algebra by defining  $v^{k}[B]=[B(k)]$, for all $k\in \Z$. The following result relates $[\mathbb{S}\mbox{Bim}]$ with $\mathcal{H}$. It is known as Soergel's categorification theorem.

\begin{teo} \cite[Theorem 1.10]{soergel2007kazhdan} \label{soergel categorificatio theorem}
There is a $\Z[v,v^{-1}]$-algebra isomorphism
\begin{equation}
\epsilon:  \mathcal{H}   \rightarrow [\mathbb{S}\mbox{Bim}]
\end{equation}
that is uniquely determined by $\epsilon(v) = R(1)$ and $\epsilon (\underline{H}_s) = B_{s}$.
\end{teo}

In order to describe an explicit inverse for $\epsilon$, we need to introduce standard bimodules. Given $x\in W$, we define the \emph{standard bimodule} $R_x$ as the $(R,R)$-bimodule such that as left $R$-module $R_{x} \cong R$, and the right action of $R$ on $R_{x}$ is determined by deforming the usual right multiplication on $R$ by $x$, i.e.,
\begin{equation} \label{right action standard bimodules}
r\cdot r' := rx(r') \mbox{ for } r\in R_x \mbox{ and } r'\in R.
\end{equation}

\begin{teo}  \cite[Theorem 5.3]{soergel2007kazhdan}   \label{inverse categorification}
The homomorphism $\epsilon: \mathcal{H} \rightarrow [\mathbb{S}Bim]$ admits an inverse, $\eta:[\mathbb{S}Bim] \rightarrow \mathcal{H} $, given by
\begin{equation}
\eta([B]) = \sum_{x\in W}  \gd (\hm^{\Z}(B, R_x)\otimes_{R}\R)   H_{x},
\end{equation}
where $\gd (-)$ signifies the respective graded dimension.
\end{teo}

Theorem \ref{inverse categorification} implies that when we expand $\eta ([B])$ in terms of the standard basis $\{ H_{w} \}_{w\in W}$ of $\mathcal{H}$, the Laurent polynomials appearing in this expansion contain positive coefficients. The following result is known in literature as \emph{Soergel's Conjecture}.

\begin{teo}  \cite[Theorem 1.1]{elias2012hodge}   \label{soergel conjecture}
For any Coxeter system $(W,S)$, we have
\begin{equation}
\epsilon (\sub{H}_{w})= [B_{w}].
\end{equation}
Consequently, $h_{x,y}= \gd (\hm^{\Z}(B_{y}, R_x)\otimes_{R}\R)$, and the KL polynomials $h_{x,y}$ have positive coefficients.
\end{teo}

\section{Libedinsky double leaves}

Let $\sub{w}$ and $\sub{v}$ be two expressions. In this section, we recall the construction of the double leaves basis (DLB) presented in \cite{libedinsky2015light}, which is a basis for $\mbox{Hom}^{\Z} (BS(\sub{w}),BS(\sub{v}))$. An alternative description of the DLB, in terms of diagrams, can be found in \cite{elias2013soergel}. It should be noted that, in our construction of the DLB, we work from right to left rather than applying the typical left-to-right approach.

\medskip
Suppose $\alpha_{s} \in V^\ast$ is an equation that belongs to a hyperplane fixed by $s \in S$. Then, for all $s \in S$, we define the \emph{Demazure operator} $\partial_{s}: R (2) \rightarrow R^s$ to be the morphism of graded $R^{s}$-modules given by
\begin{equation}
\partial_{s}(f) = \frac{f - s \cdot f}{\alpha_{s}}.
\end{equation}

In order to introduce the DLB, we first define a family of basic morphisms between Bott-Samelson bimodules. These morphisms are listed in Table \ref{table basic morphisms}, and their names, formulas, and degrees are provided. All elements belonging to the DLB can be constructed by compounding and tensoring the morphisms in this list. Therefore, the morphisms in Table \ref{table basic morphisms} generate all the morphisms that exist between Bott-Samelson bimodules. This result is due to  Libedinsky \cite[Theorem 5.1]{libedinsky2008categorie}.

\begin{table}[h]
\begin{center}
{\renewcommand{\arraystretch}{1.4}
\begin{tabular}{|l|l|c|}
  \hline
Name  & Formula & Degree \\     \hline
 $m_s: B_s \rightarrow R$   &  $1\otimes 1 \rightarrow 1 $ &   $1$ \\  \hline
  $ e_{s}:  R \rightarrow B_s  $  &  $1 \rightarrow \frac{1}{2}( \alpha_s \otimes 1 + 1 \otimes \alpha_s)$ &  $1$  \\    \hline
 $ j_s:B_s\otimes_{R} B_s \rightarrow B_s $   & $1\otimes f \otimes 1 \rightarrow \partial_{s} (f) \otimes 1$ &  $-1$\\ \hline
  $ p_s :  B_s \rightarrow B_s \otimes B_s $  &  $1 \otimes 1 \rightarrow 1\otimes 1 \otimes 1  $ &  $-1$ \\ \hline
   $ f: R \rightarrow R$  & $1\rightarrow f$  &  $\deg (f)$ \\ \hline
  $\mathbb{I}_s: B_s \rightarrow B_s   $  &  $1\otimes 1 \rightarrow 1\otimes 1  $  & 0  \\ \hline
 $ f_{st}:  B_sB_t \ldots \rightarrow  B_tB_s \ldots   $ &  &   0 \\     \hline
\end{tabular}   }
\end{center}
\caption{Basic morphisms between Bott-Samelson bimodules}
\label{table basic morphisms}
\end{table}

In the last row of Table \ref{table basic morphisms}, $B_sB_t \ldots $ (resp. $ B_tB_s \ldots  $) represents the alternating tensor product of $B_s$ and $B_t$  (resp. $B_t$ and $B_s$) having $m_{st}$ factors. Note that the formula column is empty for this row because an explicit formula (at least, not a simple formula) does not exist for this type of morphism. A morphism of type $f_{st}$ can be uniquely determined as the only degree-zero morphism from
$B_sB_t \ldots $ to $B_tB_s \ldots $, sending $1\otimes 1 \ldots \otimes 1$ to $1\otimes 1 \ldots \otimes 1$.

\begin{rem}  \rm
Subsequently, we exploit the symbol $\mathbb{I}$ in order to denote the identity of different Bott-Samelson bimodules. However, the subindex and superindex accompanying $\mathbb{I}$, as well as the context, should help avoid confusion.
\end{rem}

For each expression $\sub{w} = s_1 \ldots s_n \in S^{n}$, we inductively define a perfect binary directed tree, denoted as $\mathbb{T}_{\sub{w}}$, with nodes colored by Bott-Samelson bimodules and edges colored by morphisms from parent nodes to child nodes. In order to convert the construction of $\mathbb{T}_{\sub{w}}$ into an algorithm, we should establish the following.

\begin{cho}  \label{choices} \rm
\begin{enumerate}
\item For each $x\in W$, we fix a reduced expression, which we write as $\fre{x}$.
\item For each reduced expression $\sub{x}$ of $x\in W$, we fix a sequence of braid moves that converts $\sub{x}$ into $\fre{x}$. Then, we denote by
 $F(\sub{x}, \fre{x}): \bs{x} \rightarrow BS(\fre{x})$  the morphism obtained by replacing each braid move in the above sequence by its corresponding
 morphism of type $f_{sr}$. In particular, we fix $F(\fre{x}, \fre{x})=\mathbb{I}_{\fre{x}}$.
\item Let $x\in W$ and $s\in S$ with $l(sx)<l(x)$. Then, we fix a reduced expression $\sub{x}_s$ of $x$, which begins in $s$. We also set a sequence of braid moves that converts $\fre{x}$ into $\sub{x}_s$. Finally, we fix a morphism $F(\fre{x}, \sub{x}_s): BS(\fre{x}) \rightarrow   BS(\sub{x}_s)$ by replacing each braid move in the above sequence by its respective morphism of type $f_{sr}$.
\end{enumerate}
\end{cho}

Let us return to the construction of tree $\mathbb{T}_{\sub{w}}$. At depth $1$, the tree is the one displayed in Figure \ref{uno}.
\psset{unit=0.5}

\begin{figure}[h]
\begin{center}
\levelonetree
\end{center}
\caption{Level one of $\mathbb{T}_{\sub{w}}$}
\label{uno}
\end{figure}

Suppose $1 < k \leq n $, and assume we have constructed the tree to level $k - 1$. Let $x\in W$ be such that a node $N$ of depth $k-1$ is colored by the Bott-Samelson bimodule $\bsa{s_{1}\ldots s_{n-k+1}}\bsa{\fre{x}}$. Then, two possibilities exist.

\begin{description}
\item[a) $l(s_{n-k+1}x) > l(x) $.] Figure \ref{dos} illustrates the construction of the child nodes and edges of $N$ for this case. Here, $\sub{s_{n-k+1}x}$ represents the reduced expression of $s_{n-k+1}x$ obtained from $\fre{x}$, in which $s_{n-k+1}$ has been positioned left to $\fre{x}$.
\begin{figure}[h]
\begin{center}
\leveltwotree
\end{center}
\caption{Level $k$ of $\mathbb{T}_{\sub{w}}$}
\label{dos}
\end{figure}

\item[b) $l(s_{n-k+1}x ) < l(x) $.] In this case, the child nodes of $N$ are colored by the two Bott-Samelson bimodules located at the bottom of Figure \ref{tres}, and the child edges are colored by the morphisms obtained by composing the arrows. In this figure, the symbol $\sub{s_{n-k+1}x}$ represents the reduced expression of $s_{n-k+1}x$ obtained from $\sub{x}_{s_{n-k+1}}$ by dropping $s_{n-k+1}$ to the left of $\sub{x}_{s_{n-k+1}}$.

\begin{figure}[h]
\begin{center}
\levelthreetree
\end{center}
\caption{Level $k$ of $\mathbb{T}_{\sub{w}}$}
\label{tres}
\end{figure}

\end{description}

This completes the construction of $\mathbb{T}_{\sub{w}}$. Examining this construction, we see that all leaves of the tree  $\mathbb{T}_{\sub{w}}$ are colored by  Bott-Samelson bimodules of form $\bsa{\fre{x}}$, where $\fre{x}$ is the reduced expression fixed in Choices \ref{choices}. Furthermore, if $\sub{w}$ is a reduced expression of some element $w\in W$, the leaves in $\mathbb{T}_{\sub{w}}$ are colored by Bott-Samelson bimodules $\bsa{\fre{x}}$ with $x\leq w$.

\medskip
Composing the corresponding arrows allows us to consider each leaf in $\mathbb{T}_{\sub{w}}$ that is colored by $\bsa{\fre{x}}$ as a morphism from $\bs{w}$ to $\bsa{\fre{x}}$. Let $\mathbb{L}_{\sub{w}}(x)$ be the set of all leaves colored by $\bsa{\fre{x}}$. As mentioned earlier, we consider the set $\mathbb{L}_{\sub{w}}(x)$ to be a subset of $\hm^{\Z}(\bs{w}, \bsa{\fre{x} })$. Note that every leaf is a homogeneous morphism because all leaves were constructed as compositions of homogeneous morphisms. In fact, the degree of each leaf can be computed using $+1$ or $-1$ for each occurrence of a morphism of type $m_s$ or $j_s$, respectively.

\begin{rem}  \label{dependence of the choice} \rm
 The set $\mathbb{L}_{\sub{w}}(x)$ is not uniquely determined because it relies heavily on Choices \ref{choices}. In other words, different choices produce different morphisms in $\mathbb{L}_{\sub{w}}(x)$. Thus, when we refer to this set, we must remember that we are considering a particular fixed choice for all the non-canonical steps in the construction of $\mathbb{T}_{\sub{w}}$. For each acceptable construction of $\mathbb{T}_{\sub{w}}$, the leaves of this tree can be used to build a double leaves basis. However, different choices produce different double leaves bases.
\end{rem}

\medskip
In order to complete the construction of the double leaves basis, we need to define the adjoint tree $\mathbb{T}_{\sub{w}}^a$ of $\mathbb{T}_{\sub{w}}$. First, we note that only morphisms of type $m_s$, $j_s$, and $f_{sr}$ were used in order to construct $\mathbb{T}_{\sub{w}}$. Given any leaf $l:\bs{w} \rightarrow \bsa{\fre{x}}$ in $\mathbb{T}_{\sub{w}}$, we define its adjoint leaf  $l^a:\bsa{\fre{x}} \rightarrow \bs{w}$ by replacing each morphism  $m_s$, $j_{s}$, and $f_{sr} $ in the construction of $l$ by $e_s$, $p_s$, and $f_{rs}$, respectively. Thus, we obtain an inverted tree, $\mathbb{T}_{\sub{w}}^{a}$, which has the same nodes as $\mathbb{T}_{\sub{w}}$ but arrows pointing in opposite directions.

\medskip
Suppose $f \in \hm^{\Z}(\bs{w},\bs{u} )$ and $g\in \hm^{\Z}(\bs{x}, \bs{y})$. Then, we define
\begin{equation}
f\cdot g =\left\{
            \begin{array}{ll}
             f \circ g  , & \mbox{ if } \sub{w} = \sub{y}; \\
             \emptyset , & \mbox{ if } \sub{w} \neq \sub{y}.
            \end{array}
          \right.
\end{equation}

Given an expression $\sub{w}$, we use $\mathbb{L}_{\sub{w}}$ (resp. $\mathbb{L}_{\sub{w}}^{a}$) to indicate the set of all leaves in $\mathbb{T}_{\sub{w}}$ (resp. $\mathbb{T}_{\sub{w}}^a$). Finally, we are in a position to define the double leaves basis.

\begin{teo} \cite[Theorem 3.2]{libedinsky2015light}
For all expressions $\sub{w}$ and $ \sub{v}$, the set $\mathbb{L}_{\sub{v}}^a  \cdot \mathbb{L}_{\sub{w}}$ is a basis of $\hm^{\Z}(\bs{w}, \bs{v})$ as right $R$-module. We call this set the double leaves basis (DLB).
\end{teo}

Aside from their usefulness in defining the DLB, the leaves also describe a basis for spaces of type $\hm^{\Z}(\bs{w},R_{x})$. Note that spaces of this form appear in Theorem \ref{inverse categorification}. These spaces are relevant because they provide an alternative description for the cell modules defined in the next section. In order to explore this further, however, we must first introduce a new morphism.

For all $s \in S$, consider the $(R, R)$-bimodule morphism $\beta_s: B_s\rightarrow R_s$, determined by $\beta_s (p \otimes q) = ps(q)$, for all $p, q \in R$. On the other hand, we have $R_xR_y \cong R_{xy}$, for  all $x, y \in W$. Therefore, we can also define a morphism, denoted by $\beta_{\sub{x}} :\bs{x} \rightarrow R_{x}$, for all $x\in W$. Furthermore, let us define the set
\begin{equation}  \label{standard basis equation}
\mathbb{L}_{\sub{w}}^{\beta}(x) = \{  \beta_{\sub{x}} \circ l \mbox{ }| \mbox{ } l \in \mathbb{L}_{\sub{w}}(x)  \} \subset \hm^{\Z} (\bs{w}, R_x).
\end{equation}
In keeping with Libedinsky's work, we call $\mathbb{L}_{\sub{w}}^{\beta}(x)$ the standard leaves basis. The following proposition justifies this name choice.

\begin{pro} \cite[Proposition 6.1]{libedinsky2015light} \label{standard basis}
If $\sub{w}$ is an expression and $x\in W$, then $\mathbb{L}_{\sub{w}}^{\beta}(x)$ is an $R$-basis of $\hm^{\Z} (\bs{w}, R_x)$ as a right $R$-module.
\end{pro}

\begin{cor} \label{corollary standard basis}
Assume $\sub{w}$ is a reduced expression for $w$, and let $x\in W$. Then, $\hm^{\Z} (\bs{w}, R_x) \neq 0$ if and only if $x\leq w$.
\end{cor}
\begin{dem}
This result is a direct consequence of Proposition \ref{standard basis} if we note that $\mathbb{L}_{\sub{w}}(x)\neq \emptyset$ if and only if $x\leq w$.
\end{dem}

\section{Cellularity and KL polynomials}
Throughout this section, we fix the reduced expression $\sub{w}$. We recall that $\K$ can be considered as a left $R$-module via the isomorphism $\K \cong R/R^{+}$.  Define $\al{w}:= \mbox{End}^{\Z}(\bs{w})\otimes_{R} \K $. The $\R$-algebra $\al{w}$ has a graded cellular algebra structure --- a concept defined by Hu and Mathas \cite{hu2010graded}, who extended the work of Graham and Lehrer \cite{graham1996cellular} --- and the DLB as its graded cellular basis. We will not delve into the details of graded cellular algebras because concepts involving $\al{w}$ can be explicitly described. For a proof of the graded cellularity of $\al{w}$ the reader is referred to \cite[Proposition 6.22]{elias2013soergel} and \cite[Theorem 4.1]{plaza2014graded}.

\medskip
Graded cellular algebras are equipped with a family of modules, which are known as graded cell modules and have a bilinear form. The quotient of a graded cell module by the radical of its bilinear form is zero or simple, and all simples can be obtained in this way. In this section, we  define graded cell modules and graded simple modules for $\al{w}$ in terms of the DLB, as well as the graded decomposition numbers. We conclude this section by proving that the graded decomposition numbers and KL polynomials coincide, which is a generalization of \cite[Theorem 4.8]{plaza2014graded}.

\subsection{Graded cell modules and graded simple modules}

In order to provide a precise description for the graded cell modules and graded simple modules, we first need a couple of definitions.

\begin{defi} \label{factors through} \rm
Given two expressions $\sub{w}$ and $\sub{v}$, and assuming $x\in W$, we say that a double leaf $l_1^a \circ l_2 \in \mathbb{L}_{\sub{v}}^a  \cdot \mathbb{L}_{\sub{w}}$ \emph{factors through} $x$ whenever $l_1\in \mathbb{L}_{\sub{v}}(x)$ and $l_2 \in\mathbb{L}_{\sub{w}}(x) $.
\end{defi}

\begin{defi} \label{factors through idel} \rm
Suppose $\sub{w}$ and $\sub{v}$ are expressions. For $x\in W$, we define $\mathbb{DL}_{< x}(\sub{w},\sub{v})$ to be the $R$-submodule of  $\hm^{\Z}(\bs{w},\bs{v})$ generated by double leaves in $\mathbb{L}_{\sub{v}}^a  \cdot \mathbb{L}_{\sub{w}}$ that factor through $u < x$.
\end{defi}

We are now in a position to define the graded cell modules for $\al{w}$. Suppose $x\leq w$, $l\in \mathbb{L}_{\sub{w}}(x)$, and $a\in \mbox{End}^\Z(\bs{w})$. Then, $l\circ a$ is an element of $\hm^{\Z}(\bs{w},\bsa{\fre{x}})$, where $\fre{x}$ is the reduced expression for $x$ that was fixed in Choices \ref{choices}. Hence, we can write $l \circ a$ in terms of the DLB, and
\begin{equation}  \label{formula to define graded cell modules}
l\circ a \equiv \sum_{g\in \mathbb{L}_{\sub{w}}(x)} g r_{g} \quad \mbox{ mod } \mathbb{DL}_{<x}(\sub{w},\sub{x}),
\end{equation}
for some $r_{g} \in R$. Then, the \emph{graded cell module} $\cell{w}{x}$ is defined as the graded $\K$-vector space with basis $\mathbb{L}_{\sub{w}}(x)$. The grading on this basis is given by the grading on the leaves, and the $\al{w}$-action on $\cell{w}{x}$ is given by
\begin{equation}
l\cdot (a\otimes 1) =\sum_{g\in \mathbb{L}_{\sub{w}}(x)} g \hat{r_{g}},
\end{equation}
for every $l\in \mathbb{L}_{\sub{w}}(x)$ and $a \in \mbox{End}^\Z(\bs{w})$. Here, $\hat{r_g}\in \K$ denotes the reduction modulo $R^{+}$ of the scalars, $r_g$, that appear in (\ref{formula to define graded cell modules}). The following proposition provides an alternative description of the cell modules.

\begin{pro} \cite[Lemma 4.6]{plaza2014graded} \label{isomorphism cell}
Assume $\sub{w}$ is a reduced expression for $w \in W$. Then, for all $x \leq w$, we have an isomorphism
\begin{equation}
\hm^{\Z}(\bs{w}, R_x)\otimes_{R} \K \cong \cell{w}{x}
\end{equation}
of right $\al{w}$-modules, where the $\al{w}$-action on $\hm^{\Z}(\bs{w}, R_x)\otimes_{R} \K $ is given by the composition of morphisms.
\end{pro}

Now we can define a bilinear form, $\langle \mbox{ , } \rangle$, on $\cell{w}{x}$. Let $l_1, l_2 \in \mathbb{L}_{\sub{w}}(x)$. The composition $l_1 \circ l_2^{a}$ belongs to  $\mbox{End}^\Z (\bsa{\fre{x}})$, and we can expand it in terms of the DLB as
\begin{equation}
l_1 \circ l_2^{a} \equiv \mathbb{I}_{\fre{x}} r(l_1,l_2) \quad \mbox{ mod } \mathbb{DL}_{<x}(\sub{x} , \sub{x}),
\end{equation}
for some $r(l_1,l_2) \in R$. Then, we define the bilinear form on $\cell{w}{x}$ to be
\begin{equation}
\langle l_1 , l_{2} \rangle = \widehat{r(l_1,l_2)}  \in \K,
\end{equation}
where $\widehat{r(l_1,l_2)}$ denotes  the reduction modulo $R^{+}$ of $r(l_1,l_2)$. This bilinear form is symmetric, associative and homogeneous, and its radical is defined to be
\begin{equation}
\mbox{Rad} (\cell{w}{x})= \{l \in \cell{w}{x} \mbox{ }|\mbox{ } \langle l,l'\rangle =0 \mbox{ for all } l' \in \cell{w}{x} \}.
\end{equation}
An easy consequence of the associativity and homogeneity of $\langle \mbox{ , } \rangle $ is that $\mbox{Rad} (\cell{w}{x})$ is a graded $\al{w}$-submodule of $\cell{w}{x}$. Therefore, we can define the quotient $\al{w}$-module $\simple{w}{x}:= \cell{w}{x} / \mbox{Rad} (\cell{w}{x})$. The modules $\{ \simple{w}{x} \}_{x\leq w}$ are simple or zero.  The bilinear form controls the behavior of the category of Soergel bimodules, as explicitly stated in the following lemma.

\begin{lem}   \label{decomposition of Bott-Samelson bimodules}
Suppose that $\sub{w}$ is a reduced expression for $w\in W$, and define
\begin{equation}
\simpleset{w}:= \{x \leq w  \mbox{ } | \mbox{  } \simple{w}{x}\neq 0 \}.
\end{equation}
Then, there exists an isomorphism of graded $(R,R)$-bimodules:
\begin{equation}
\bs{w} = \bigoplus_{y\in \simpleset{w}} \gd \simple{w}{y} B_{y}.
\end{equation}
\end{lem}
\begin{dem}
A proof for this can be found in \cite[Lemma 4.5]{plaza2014graded}. However, we remark that this proof rephrases the proof of \cite[Lemma 3.1]{williamson2013schubert} using the language of graded cellular algebras.
\end{dem}

\medskip
Given a graded $\R$-vector space $M$ and $k \in \Z$, we let $M\langle k \rangle$ represent the graded $\K$-vector space obtained from $M$ by shifting the grading on $M$, i.e., $M\langle k \rangle_{i} = M_{i - k}$ for all $i \in \Z$.

\begin{rem} \rm
In Section 2.2,  the grading shift on graded $(R,R)$-bimodules was defined using round brackets. In the previous paragraph, we have defined the grading shift on graded $\R$-vector spaces by using angled brackets. These two concepts differ only in the direction of the grading shift. In fact, we could have used round  or angled brackets to define both concepts. However, we prefer to distinguish them to keep the notation as consistent as possible with both the literature about Soergel bimodules and graded cellular algebras.
\end{rem}

The following theorem classifies the set of simple graded $\al{w}$-modules.

\begin{teo} \cite[Theorem 3.4]{graham1996cellular}, \cite[Theorem 2.10]{hu2010graded} \label{clasification of graded simple modules}
If $\sub{w}$ is a reduced expression, then the  set $$\{ \simple{w}{x}\langle k \rangle \mbox{ } |\mbox{ } x\in \simpleset{w} \mbox{ and } k\in \Z \}$$ is a complete set of pairwise non-isomorphic graded simple right $\al{w}$-modules.
\end{teo}

\subsection{Graded decomposition numbers and KL polynomials}

Suppose $A$ is a graded algebra. If $M$ is a graded $A$-module and $L$ is a graded simple $A$-module, let $[M:L\langle k \rangle]$ be the multiplicity of the simple module $L\langle k \rangle$ as a graded composition factor of $M$, for all $k\in \Z$. The graded decomposition number is defined to be
\begin{equation}
d_{A}(M , L) = \sum_{k\in \Z}  [M : L \langle k \rangle]v^k \in \Z[v,v^{-1}].
\end{equation}
In particular, if $A=\al{w}$, $x\leq w$, and $y\in \simpleset{w}$, we will write
\begin{equation}
\gdn{w}{x}{y}:= d_{\al{w}}(\cell{w}{x},\simple{w}{y} ).
\end{equation}
We refer to $\{\gdn{w}{x}{y} \mbox{ } | \mbox{ } x \leq w \mbox{ and } y\in \simpleset{w} \} $ as the graded decomposition numbers for $\al{w}$. The general theory of graded cellular algebras tells us that the graded decomposition numbers for $\al{w}$ satisfy the following triangularity property \cite[Lemma 2.13]{hu2010graded}:
\begin{equation}  \label{triangularity property for gdn}
\gdn{w}{y}{y}= 1 \mbox{ and }  \gdn{w}{x}{y} = 0 \mbox{ unless  } x\leq y,
\end{equation}
for every $x\leq w$ and $y\in \simpleset{w}$. The following lemma helps relate the graded decomposition numbers for $\al{w}$ with the KL polynomials. All the statements in the lemma are well-known, see for example \cite[Appendix]{donkin1998q} for the ungraded case.

\begin{lem}  \label{lema gdn equals in idempotent algbera}
Assume $A$ is a graded algebra and $e\in A$ is a homogeneous idempotent. Then the idempotent subalgebra $eAe$ of $A$ is again a graded algebra, and, for any graded right $A$-module $M$, $Me$ has a natural structure of  graded right $eAe$-module. Moreover, if $\{ L_{\lambda} \mbox{ } | \mbox{ } \lambda \in \Lambda \}$ is a complete (up to a degree shift) set of representatives of non-isomorphic classes of graded simple $A$-modules, then we have the following.
\begin{enumerate}
  \item The set $\{ L_{\lambda}e \mbox{ } | \mbox{ } \lambda \in \Lambda  \mbox{ and } L_{\lambda}e\neq 0 \}$ is a complete (up to a degree shift) set of representatives of non-isomorphic classes of graded simple $eAe$-modules.
  \item If $M$ is a graded right $A$-module and $L$ is a graded right simple $A$-module such that $Le\neq 0$, then
\begin{equation}   \label{equation gdn equal in idempotent algebra}
d_{A}(M,L)=d_{eAe}(Me,Le).
\end{equation}
\item Furthermore, if $e$ is a primitive idempotent, then there exists exactly one (up to a degree shift) graded right simple $eAe$-module. This unique simple module must be of the form $Le$ for some simple graded $A$-module $L$.
\end{enumerate}
\end{lem}

\begin{pro}    \label{equality graded decompostion numbers with KL}
Let $\sub{w}$ be a reduced expression for $w\in W$. Then
\begin{equation}
 \gdn{w}{x}{y} = h_{x,y},
\end{equation}
for all $x\leq w$ and $y \in \simpleset{w}$.
\end{pro}
\begin{dem}
First, we recall from Theorem \ref{soergel conjecture} that $h_{x,y}=\gd (\hm^{\Z} (B_{y}, R_{x}) \otimes_{R} \K)$. In addition, Lemma \ref{decomposition of Bott-Samelson bimodules} states that, for all $y\in \simpleset{w}$, we can choose a primitive idempotent, $\epsilon_{y}^{k} \in \mbox{End}^{\Z}(\bs{w})$, whose image is isomorphic to $B_{y}(k)$ for some $k\in \Z$. Reduction modulo $R^{+}$ of $\epsilon_{y}^{k}$, which we simply denote using $e$, is again a primitive idempotent in $\al{w}$. According to Lemma \ref{lema gdn equals in idempotent algbera}, the idempotent subalgebra $e \al{w} e $ is a graded algebra with a unique (up to a degree shift) simple graded right module. Furthermore, by Theorem \ref{clasification of graded simple modules} and Lemma \ref{lema gdn equals in idempotent algbera}, we know that the unique simple graded $e\al{w}e$-module is of the form $\simple{w}{z} e \neq 0$, for some $z\in \simpleset{w}$. We claim that $z=y$. To prove this claim, for $x\leq w$, we have
\begin{equation}  \label{list of equations}
\begin{array}{rl}
 v^{k} \gd (\hm^{\Z}(B_y,R_x) \otimes_R \K)  &= \gd (\hm^{\Z}(B_y(k),R_x) \otimes_R \K) \\
   & =\gd \left( (\hm^{\Z}(\bs{w},R_x) \otimes_R \K) e \right)\\
   & =\gd (\cell{w}{x} e)  \\
   & =d_{e\al{w}e}(\cell{w}{x} e, \simple{w}{z} e  ) \gd ( \simple{w}{z}  e) \\
   & =\gdn{w}{x}{z} \gd (\simple{w}{z} e).
\end{array}
\end{equation}
Here, the second equation is a result of the definition of $e$, the third equation is a consequence of Proposition \ref{isomorphism cell}, the fourth equation follows from the fact that $\simple{w}{z} e$ is the unique simple $ e \al{w} e$-module, and the last equation is obtained using Lemma \ref{lema gdn equals in idempotent algbera}. Substituting $x=y$ into (\ref{list of equations}) gives us
\begin{equation}   \label{equation intermediate}
v^{k} \gd (\hm^{\Z}(B_y,R_y) \otimes_R \K) = \gdn{w}{y}{z} \gd (\simple{w}{z} e).
\end{equation}
Because $\gd (\hm^{\Z}(B_y,R_y) \otimes_R \K)=h_{y,y}=1 $, the right side of (\ref{equation intermediate}) is not zero. Therefore, $\gdn{w}{y}{z}\neq 0$, and applying (\ref{triangularity property for gdn}) gives us $y\leq z$. On the other hand, substituting $x=z$ into (\ref{list of equations}) provides
\begin{equation}  \label{equation intermediate two}
v^{k} \gd (\hm^{\Z}(B_y,R_z) \otimes_R \K) = \gdn{w}{z}{z} \gd (\simple{w}{z} e).
\end{equation}
Again, by (\ref{triangularity property for gdn}) we know that $\gdn{w}{z}{z}=1$. Thus, we conclude that $h_{z,y}=\gd (\hm^{\Z}(B_y,R_z) \otimes_R \K)\neq 0$, and $y\geq z$. This proves our claim.

Furthermore, if we substitute $y=z$ into (\ref{list of equations}), we obtain
\begin{equation} \label{equation intermediate three}
v^{k} \gd (\hm^{\Z}(B_y,R_x) \otimes_R \K) = \gdn{w}{y}{x} \gd (\simple{w}{y} e).
\end{equation}
Finally, substituting $x=y$ into (\ref{equation intermediate three}) provides $\gd (\simple{w}{y} e)=v^{k}$, and cancelling these terms in $(\ref{equation intermediate three})$ produces the desired equality.
\end{dem}

\section{Branching rules}

In order to proceed, we fix from now on a reduced expression $\sub{w}=s_{1} \ldots s_{k}$ of $w\in W$. In addition, in order to simplify the notation, we refer to $s_1$ simply as $s$. Define $x'=sx \in W$, for all $x\in W$. In particular, we have $w'=sw$. From now on we fix the reduced expression $\sub{w'}=s_{2} \ldots s_{k}$ for $w'$. By definition of $\al{w}$ and $\al{w'}$ it is clear that there exists an injective $\R$-algebra homomorphism from $\al{w'}$ to $\al{w}$, which sends an element $a\in \al{w'} $ to $ \mathbb{I}_{s} \otimes a \in \al{w}$. This embedding allows us to view $\al{w'}$ as a subalgebra of $\al{w}$. Accordingly, any $\al{w}$-module $M$ can be considered as an $\al{w'}$-module by restriction. We denote the restricted module by $\res{M}$. In this section, we obtain graded branching rules for graded cell modules and graded simple modules of $\al{w}$.  We end this section by showing how to utilize these graded branching rules to recover a normalized version of (\ref{equation introduction}).

\subsection{Branching rules for cell modules}

Consider the expression $\sub{w'}$ and its corresponding tree $\mathbb{T}_{\sub{w'}}$. Recall that in order to begin the construction of $\mathbb{T}_{\sub{w'}}$ we must fix some data (see Choices \ref{choices}). In particular, the reduced expression $\fre{x}$ must be fixed, for all $x \in W$. There is no loss of generality in assuming that $s$ appears on the left side of $\fre{x}$, for every $x$ with $x'<x$. After the reduced expression has been fixed for every element $x$ with $x'<x$, we can fix reduced expressions $\fre{y}$ for every element $y\in W$, where $y'>y$, according to the following rule:  $\fre{y}$ is obtained from $\fre{y'}$ by deleting the $s$ located to the left of $\fre{y'}$. It should be noted that $\fre{y'}$ has already been fixed since $sy'<y'$.

\begin{figure}[h]
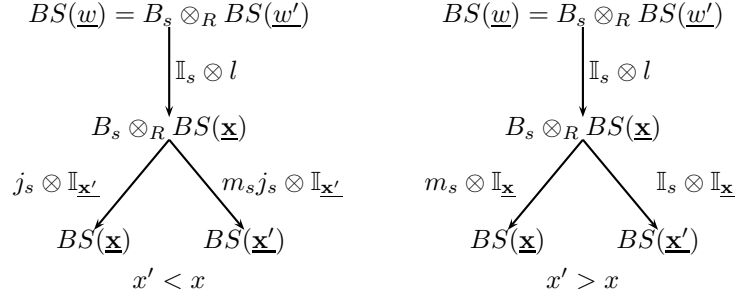

\begin{center}
\newleavesfromoldleaves
\end{center}
\caption{Constructing leaves from $\mathbb{T}_{\sub{w'}}$ to $\mathbb{T}_{\sub{w}}$}
\label{building leaves}
\end{figure}

\medskip
The consequence of these choices is that the tree $\mathbb{T}_{\sub{w}}$ can be obtained directly from $\mathbb{T}_{\sub{w'}}$. Concretely, each leaf $l\in \mathbb{L}_{\sub{w'}}(x)$ gives rise to two leaves in $\mathbb{T}_{\sub{w}}$, according to $x'<x$ or $x'>x$, which is illustrated in Figure \ref{building leaves}. Although this could be considered a direct consequence of the construction of $\mathbb{T}_{\sub{w}}$, we note the absence of morphisms of type  $F(\mbox{ , } )$ in Figure \ref{building leaves}, which results from the previously realized choices.

\begin{lem}  \label{partitionating leaves}
Given $x\leq w$, the set $\mathbb{L}_{\sub{w}}(x)$ can be partitioned into two disjoint sets, $ Y_{\sub{w},1}(x)$ and $ Y_{\sub{w},2}(x)$, which are defined as follows.
\begin{equation}
\begin{array}{l}
 Y_{\sub{w},1}(x)=\left\{
     \begin{array}{ll}
       \left\{ (\mathbb{I}_{s} \otimes l) \mbox{ } | \mbox{ } l\in \mathbb{L}_{\sub{w'}}(x') \right\}, & \mbox{if  } x'<x;\vspace{.3cm} \\
       \left\{ (m_s \otimes l)    \mbox{ } | \mbox{ } l\in \mathbb{L}_{\sub{w'}}(x) \right\}   , &  \mbox{if } x'>x
     \end{array}
   \right.
\\

\\
 Y_{\sub{w},2}(x)= \left\{
     \begin{array}{ll}
      \left\{ (j_s \otimes \mathbb{I}_{\fre{x'}})(\mathbb{I}_{s} \otimes l)\mbox{ }|\mbox{ }l\in \mathbb{L}_{\sub{w'}}(x) \right\} , &\mbox{if  } x'<x; \vspace{.3cm} \\
  \left\{ (m_s j_s \otimes \mathbb{I}_{\fre{x'}})(\mathbb{I}_{s} \otimes l) \mbox{ } | \mbox{ } l\in \mathbb{L}_{\sub{w'}}(x') \right\} , & \mbox{if } x'>x
     \end{array}
   \right.
\end{array}
\end{equation}
\end{lem}
\begin{dem}
The reason for this is clear given the construction of the leaves in $\mathbb{T}_{\sub{w}}$ from the leaves in $\mathbb{T}_{\sub{w'}}$, which is provided in Figure \ref{building leaves}.
\end{dem}

\medskip
The branching rules for $\cell{w}{x}$ can be classified into two cases: $x'<x$ and $x'>x$. In the following theorem, we assume $x'<x$.

\begin{teo} \label{teo short exact sequence x greater than xs}
Suppose $x\leq w$ is such that $x'<x$. Then, there exists a short exact sequence of graded $\al{w'}$-modules,
\begin{equation}\label{exact sequence x greater than xs}
0 \longrightarrow \hm^{\Z}(\bs{w' },R_{x'} ) \otimes_{R} \K \overlim{\phi_{1}}   \res{\hm^{\Z}(\bs{w },R_{x} ) \otimes_{R} \K} \overlim{\phi_{2}} \left(\hm^{\Z}(\bs{w'}, R_{x} ) \otimes_{R} \K\right)\langle -1 \rangle \longrightarrow 0,
\end{equation}
in which $\phi_{1}(f)=\beta_{s} \otimes f$ and $\phi_{2}(f)=f (e_{s}\otimes \mathbb{I}_{\sub{w'}})$.
\end{teo}
\begin{dem}
Assume $a\in \al{w'}$. Then, for all $f\in\hm^{\Z}(\bs{w' },R_{x'} ) \otimes_{R} \K $, we have
\begin{equation}
\phi_{1}(fa)=\beta_{s}\otimes (f a) = (\beta_{s} \otimes f)(\mathbb{I}_{s} \otimes a)=\phi_{1}(f)(\mathbb{I}_{s} \otimes a),
\end{equation}
proving that $\phi_{1}$ is an $\al{w'}$-algebra homomorphism. Similarly, we have
\begin{equation}
\phi_2(f (\mathbb{I}_s \otimes a))= f(\mathbb{I}_{s} \otimes a ) (e_s\otimes \mathbb{I}_{\sub{w'}})
                                        = f  (e_s \otimes \mathbb{I}_{\sub{w'}})  a
                                        = \phi_{2}(f)a,
\end{equation}
for every $f\in \hm^{\Z}(\bs{w },R_{x} ) \otimes_{R} \K$. Hence, $\phi_{2}$ is also an $\al{w'}$-algebra homomorphism. Furthermore, we note that $\phi_{1}$ preserves the grading, and $\phi_2$ increases the grading by $1$ because $\deg (e_s)=1$. In summary, both $\phi_{1}$ and $\phi_{2}$ are homogeneous $\al{w'}$-algebra homomorphisms of degree $0$ and $1$, respectively.

\medskip
Now, we are going to prove  that $\phi_1$ and $\phi_2$ are injective and surjective, respectively. Proposition \ref{standard basis} ensures that $\mathbb{L}_{\sub{w}}^{\beta}(x)$ is an $R$-basis of $\hm^{\Z}(\bs{w },R_{x} )$. Therefore, the image of $\mathbb{L}_{\sub{w}}^{\beta}(x)$ under reduction modulo $R^+$ is an $\K$-basis of $\hm^{\Z}(\bs{w },R_{x} ) \otimes_{R} \K$. If we relax the notation, we can also consider $\mathbb{L}_{\sub{w}}^{\beta}(x)$ to be an $\K$-basis of  $\hm^{\Z}(\bs{w },R_{x} ) \otimes_{R} \K$. Lemma \ref{partitionating leaves} implies that, given $x'<x$,  $\mathbb{L}_{\sub{w}}^{\beta}(x)$ can be partitioned into two disjoint sets:
$$
\begin{array}{rl}
 Y_{\sub{w},1}^{\beta}(x):=  & \{  \beta_{\sub{\mathbf{x}}} \circ (\mathbb{I}_s \otimes l) \mbox{ } | \mbox{ } l \in \mathbb{L}_{\sub{w'}}(x')  \},  \mbox{ and }\vspace{.2cm}\\
 Y_{\sub{w},2}^{\beta}(x):=  & \{  \beta_{\sub{\mathbf{x}}} \circ  (j_s \otimes \mathbb{I}_{\sub{\mathbf{x'}}})(\mathbb{I}_s\otimes l)  \mbox{ } | \mbox{ } l \in\mathbb{L}_{\sub{w'}}(x)    \}.
\end{array}
$$

 For all $l\in \mathbb{L}_{\sub{w'}}(x')$, we have
$$
\begin{array}{rl}
  \phi_1 (\beta_{\sub{\mathbf{x'}}} \circ l) & = \beta_s \otimes (\beta_{x'} \circ l) \\
   & = (\beta_s \otimes \beta_{\sub{\mathbf{x'}}}) ( \mathbb{I}_s \otimes l) \\
   & = \beta_{\sub{\mathbf{x}}} \circ ( \mathbb{I}_s \otimes l).
\end{array}
$$
Therefore, $ \phi_1 (\mathbb{L}_{\sub{w'}}^{\beta}(x'))=  Y_{\sub{w},1}^{\beta}(x)$, and $ \phi_1 $ is injective. Similarly, for $ l\in  \mathbb{L}_{\sub{w'}}(x)$ we have

$$
\begin{array}{rl}
 \phi_2 ( \beta_{\sub{\mathbf{x}}} \circ  (j_s \otimes \mathbb{I}_{\sub{\mathbf{x'}}})(\mathbb{I}_s\otimes l))& = \beta_{\sub{\mathbf{x}}} \circ  (j_s \otimes \mathbb{I}_{\sub{\mathbf{x'}}})(\mathbb{I}_s\otimes l) (e_s \otimes \mathbb{I}_{\sub{w'}})  \\
   &    = \beta_{\sub{\mathbf{x}}} \circ  (j_s \otimes \mathbb{I}_{\sub{\mathbf{x'}}})(e_s\otimes l)  \\
   &   =   \beta_{\sub{\mathbf{x}}} \circ l
\end{array}
$$
because $j_s(e_s \otimes \mathbb{I}_{s}) = \mathbb{I}_{s}$. Hence, $\phi_2(Y_{\sub{w},2}^{\beta}(x))= \mathbb{L}_{\sub{w'}}^{\beta}(x)$, and $\phi_2$ is surjective.

\medskip
Finally, in order to complete the proof we need to show that $\mbox{Im} (\phi_1) = \ker (\phi_2) $. For each $l\in \mathbb{L}_{\sub{w'}}(x')$, we have
\begin{equation}   \label{image into the kernel}
\phi_2 \phi_1 (\beta_{\sub{\mathbf{x'}}} \circ l )  =\phi_2(\beta_{\sub{\mathbf{x}}} \circ (\mathbb{I}_{s}\otimes l ))= \beta_{\sub{\mathbf{x}}} (\mathbb{I}_{s}\otimes l )(e_s \otimes \mathbb{I}_{\sub{w'}})=
\beta_{\sub{\mathbf{x}}} (e_s \otimes \mathbb{I}_{\sub{\mathbf{x'}}} )l =0,
\end{equation}
because $\beta_{\sub{\mathbf{x}}} (e_s \otimes \mathbb{I}_{\sub{\mathbf{x'}}} ) \in \hm^{\Z}(\bs{\mathbf{x'}},R_{x} ) \otimes_{R} \K =0$, according to Corollary \ref{corollary standard basis}. Equation (\ref{image into the kernel}) implies that $\mbox{Im}  (\phi_1) \subset \ker (\phi_2) $. Moreover, an easy counting argument reveals that $\mbox{Im} (\phi_1) = \ker (\phi_2) $, completing the proof.
\end{dem}

\begin{cor} \label{cell teo short exact sequence x greater than xs}
Assume $x\leq w$ is such that $x'<x$. Then there exists a short exact sequence of graded $\al{w'}$-modules
\begin{equation}\label{cell exact sequence x greater than xs}
0 \longrightarrow \cell{w'}{x'} \longrightarrow  \res{\cell{w}{x}} \longrightarrow \cell{w'}{x}\langle -1 \rangle \longrightarrow 0.
\end{equation}
\end{cor}
\begin{dem}
This is a direct consequence of Proposition \ref{isomorphism cell} and Theorem \ref{teo short exact sequence x greater than xs}.
\end{dem}

Analogous results for Theorem \ref{teo short exact sequence x greater than xs} and Corollary \ref{cell teo short exact sequence x greater than xs} also exist for the case in which $x' > x$.

\begin{teo}  \label{teo branching rules x' smaller than x}
If $x\leq w$ is given such that $x'>x$, then there exists a graded short exact sequence of $\al{w'}$-modules,
\begin{equation}\label{exact sequence x smaller than xs}
0 \longrightarrow \left( \mbox{\small Hom}^{\Z}(\bs{w' },R_{x} ) \otimes_{R} \K \right) \langle 1 \rangle  \overlim{\phi_{1}}
\res{\hm^{\Z}(\bs{w },R_{x} ) \otimes_{R} \K} \overlim{\phi_{2}}   \hm^{\Z}(\bs{w'}, R_{x'} ) \otimes_{R} \K  \longrightarrow 0,
\end{equation}
where $\phi_{1}(f)=m_{s} \otimes f$ and $\phi_{2}(f)= ( \beta_{s} \otimes f) (p_se_s  \otimes \mathbb{I}_{\sub{w'}} ) $. Equivalently, there exists a graded short exact sequence of $\al{w'}$-modules:
\begin{equation}\label{cell exact sequence x smaller than xs}
0 \longrightarrow \cell{w'}{x} \langle 1 \rangle  \longrightarrow  \res{\cell{w}{x}} \longrightarrow \cell{w'}{x'} \longrightarrow 0.
\end{equation}
\end{teo}

\begin{dem}
This proof is analogous to that of Theorem \ref{teo short exact sequence x greater than xs}.
\end{dem}

\medskip
In order to continue, it would be useful to rephrase the branching rules in terms of Grothendieck groups. Given a reduced expression $\sub{w}$, the \emph{Grothendieck group} $\mathcal{G}_{\sub{w}}$ of $\al{w}$ is the $\Z[v,v^{-1}]$-module generated by symbols $[M]$, in which $M$ runs over all isomorphism classes of  graded $\al{w}$-modules, together with the following relations:
\begin{enumerate}
\item   $[M\langle k \rangle]= v^{k} [M]$, for every graded right $\al{w}$-module $M$ and every $k\in \Z$;
\item   $[M]=[M']+[M'']$, if there exists a (graded) short exact sequence $   0\rightarrow M' \rightarrow M \rightarrow M'' \rightarrow 0  $ of  graded right $\al{w}$-modules.
\end{enumerate}
Thus, Theorem \ref{clasification of graded simple modules} reveals that  $\mathcal{G}_{\sub{w}}$ is free as a $\Z[v,v^{-1}]$-module with a basis given by
$\{ [\simple{w}{x}] \mbox{ } | \mbox{ } x\in \simpleset{w} \}$. Using this terminology, the graded short exact sequences in Theorem \ref{teo short exact sequence x greater than xs} and Theorem \ref{teo branching rules x' smaller than x} can be rephrased as follows.

\begin{teo}  \label{teo grothendieck branching rules x greater than xs}
 There exists a homomorphism of $\Z[v,v^{-1}]$-modules $\mbox{Res}: \mathcal{G}_{\sub{w}} \rightarrow \mathcal{G}_{\sub{w'}}$  determined by
\begin{equation}
 \res{[M]} = [\res{M}],
\end{equation}
for every right graded $\al{w}$-module $M$. The image of the class of cell module $\cell{w}{x}$ under this homomorphism is
\begin{equation}  \label{grothendieck group exact sequence}
     \mbox{Res}([\cell{w}{x}] )=  \left\{
                                  \begin{array}{rl}
                                    v^{-1} [\cell{w'}{x}] + [\cell{w'}{x'}], & \mbox{if } x'<x ; \\
                                     v [\cell{w'}{x}] + [\cell{w'}{x'}], & \mbox{if }  x'>x.
                                  \end{array}
                                \right.
\end{equation}
\end{teo}

\subsection{Branching rules for simple modules}

In this section, we obtain graded branching rules for simple $\al{w}$-modules. In other terms, we calculate the image of the class $[\simple{w}{x}]$ under the action of the homomorphism $\mbox{Res}$, for all $x\in \simpleset{w}$. In order to do this, we define a $\Z[v,v^{-1}]$-homomorphism $\rho$ from $\mathcal{G}_{\sub{w}}$ to $\mathcal{G}_{\sub{w'}}$ and prove that $\rho$ coincides with $\mbox{Res}$. Given $u\in \simpleset{w'}$, there exist polynomials $h_{s,u}^{x}\in \Z[v,v^{-1}]$ such that
\begin{equation}  \label{decomposition structure constant}
B_s \otimes B_{u} \cong \bigoplus_{x\in \simpleset{w}}  h_{s,u}^{x} B_x.
\end{equation}
In principle, the sum on the right side of (\ref{decomposition structure constant}) should run over $W$. However, Lemma \ref{decomposition of Bott-Samelson bimodules} guarantees that if  $x \not \in \simpleset{w}$, then the polynomials $h_{s,u}^{x}$ equal zero.

\medskip
Using the polynomials $h_{s,u}^{x}$, we can define a  $\Z[v,v^{-1}]$-homomorphism $\rho: \mathcal{G}_{\sub{w}} \rightarrow \mathcal{G}_{\sub{w'}}$ determined in the basis $\{ [\simple{w}{x}] \mbox{ } | \mbox{ } x\in \simpleset{w} \}$ of $\mathcal{G}_{\sub{w}}$ by
\begin{equation}  \label{definition of rho}
\rho([\simple{w}{x}])= \sum_{u\in \simpleset{w'}} h_{s,u}^{x} [\simple{w'}{u}],
\end{equation}
for every $x\in \simpleset{w}$. In order to prove that $\rho = \mbox{Res}$, we need the following lemmas.

\begin{lem}   \label{preparation}
Given any $z\leq w$ and any $u\in \simpleset{w'}$, we have
\begin{equation}
d_{\al{w'}}(\mbox{Res} (\cell{w}{z}), \simple{w'}{u})= \sum_{ x\in \simpleset{w} } h_{s,u}^{x}\gdn{w}{z}{x}.
\end{equation}
\end{lem}
\begin{dem}
Assuming $u\in \simpleset{w'}$, we can choose a primitive idempotent, $\epsilon_{u}^{k} \in \mbox{End}^{\Z}(\bs{w'})$, whose image is isomorphic to $B_{u}(k)$, for some $k\in \Z$. Reduction modulo $R^{+}$ of $\epsilon_{u}^{k}$, which we denote as $e \in \al{w'}$, is a primitive idempotent, and from  Lemma \ref{lema gdn equals in idempotent algbera} we know the  algebra $e\al{w'}e$ has a unique (up to a degree shift) graded simple module. In the proof of Proposition \ref{equality graded decompostion numbers with KL}, we have already proven that the unique simple module of  $e\al{w'}e$  is $\simple{w'}{u}e$. Furthermore, we can also prove that $\gd{\simple{w'}{u}e}= v^k$.  Thus, we have the following chain of equations.

\begin{eqnarray}
 d_{\al{w'}}(\mbox{Res}(\cell{w}{z}),\simple{w'}{u}) v^{k}  &=& d_{\al{w'}}(\mbox{Res}(\cell{w}{z}),\simple{w'}{u})\gd \simple{w'}{u}e  \label{chain one}  \\
   &=& d_{e\al{w'}e}(\mbox{Res}(\cell{w}{z}) (\mathbb{I}_{s}\otimes e), \simple{w'}{u} e)\gd \simple{w'}{u} e \label{chain two} \\
   &=& \gd \cell{w}{z} (\mathbb{I}_s \otimes e) \label{chain three} \\
   &=& \gd \left( \hm^{\Z}(\bs{w},R_z)\otimes_R \K \right)(\mathbb{I}_s \otimes e) \label{chain four}  \\
   &=& v^k \gd \left( \hm^{\Z}(B_s\otimes B_{u},R_z)\otimes_R \K \right) \label{chain six} \\
   &=& v^{k} \sum_{x\in \simpleset{w}} h_{s,u}^x \gd \left( \hm^{\Z}(B_x,R_z) \otimes_R \K \right)  \label{chain seven} \\
   &=& v^{k} \sum_{ x\in \simpleset{w} } h_{s,u}^x \gdn{w}{z}{x}.  \label{chain eight}
\end{eqnarray}
Here, we obtain (\ref{chain one}) because: $\gd{\simple{w'}{u}e}= v^k$; (\ref{chain two}) follows from Lemma \ref{lema gdn equals in idempotent algbera}; (\ref{chain three}) applies the fact that $\simple{w'}{u}e$ is the unique (up to a degree shift) simple $e\al{w'}e$-module; Proposition \ref{isomorphism cell} implies (\ref{chain four}); the definition of $e$ justifies (\ref{chain six}); (\ref{chain seven}) is a consequence of (\ref{decomposition structure constant}); and (\ref{chain eight}) follows from Theorem \ref{soergel conjecture} and Proposition \ref{equality graded decompostion numbers with KL}. Finally, multiplying the chain of equations by $v^{-k}$, we achieve the desired equality.
\end{dem}

\begin{lem}   \label{graded module write in simple basis}
Given a graded right $\al{w}$-module $M$, we have the following equation in $\mathcal{G}_{\sub{w}}$:
\begin{equation}
[M]= \sum_{x\in\simpleset{w}} d_{\al{w}}(M,\simple{w}{x}) [\simple{w}{x}].
\end{equation}
Likewise, the same is true if we replace $\sub{w}$ with $\sub{w'}$.
\end{lem}
\begin{dem}
This is a direct consequence of the definitions of graded decomposition numbers and Grothendieck groups.
\end{dem}

\begin{teo}
The homomorphisms $\mbox{Res}$ and $\rho$ coincide.
\end{teo}
\begin{dem}
First, according to the triangularity property in (\ref{triangularity property for gdn}), we know that $\{ [ \cell{w}{x} ] \}_{x\in \simpleset{w}}$ is also a $\Z[v,v^{-1}]$-basis of $\mathcal{G}_{\sub{w}}$. Therefore, in order to prove the theorem, it is sufficient to show that
\begin{equation}  \label{equation to prove rho and res}
 \mbox{Res}\left([\cell{w}{x}]\right)= \rho \left([\cell{w}{x}]\right),
\end{equation}
for all $x\in \simpleset{w}$. More generally, we will prove that (\ref{equation to prove rho and res}) holds for all $x\leq w$. According to Lemma \ref{graded module write in simple basis}, given $x\leq w$, we have
\begin{equation}  \label{cell in terms of simples}
[\cell{w}{x}] = \sum_{y\in\simpleset{w}} \gdn{w}{x}{y} [\simple{w}{y}].
\end{equation}
Thus, applying  the map $\rho$ to (\ref{cell in terms of simples}) and using the definition of $\rho $ provided in (\ref{definition of rho}) gives us
\begin{equation}  \label{image of rho}
\rho \left([\cell{w}{x}]\right)=\sum_{u\in \simpleset{w'}}\left( \sum_{ y\in \simpleset{w} } h_{s,u}^y \gdn{w}{x}{y}  \right)[\simple{w'}{u}],
\end{equation}
for every $x\leq w$. On the other hand, the combination of Lemmas \ref{preparation} and \ref{graded module write in simple basis}  implies that the right side of (\ref{image of rho}) equals $\mbox{Res} ( [\cell{w}{x}])$, for all $x\leq w$. This completes the proof.
\end{dem}

\begin{cor}   \label{corollary branching simples}
If $x\in \simpleset{w}$, then the restriction of the simple module $\simple{w}{x}$ is given by
\begin{equation}
\mbox{Res}([\simple{w}{x}]) = \sum_{u\in \simpleset{w'}} h_{s,u}^{x} [\simple{w'}{u}].
\end{equation}
\end{cor}

\subsection{Categorifying a recursive formula for KL polynomials}

In this subsection, we use the branching rules for graded cell modules  and graded simple modules in order to obtain (\ref{alg kl poly cinco}). In other words, we use the categories of graded finite-dimensional $\al{w}$-modules and $\al{w'}$-modules (and, therefore, the category of Soergel bimodules) in order to recover (\ref{alg kl poly cinco}). In short, we categorify (\ref{alg kl poly cinco}). Throughout this subsection, we maintain the notation provided in this section's introduction.

In order to categorify (\ref{alg kl poly cinco}), we need the following lemma. First, however, we recall that $\mu(z,w')$ denotes the linear coefficient of the Kazhdan-Lusztig polynomial $h_{z,w'}$, for all $z\in W$ (see section \ref{hecke and kl poly}).

\begin{lem} \label{lemma to categorify}
For all $z<w$, we have
\begin{equation}  \label{equation to categorify}
h_{s,w'}^{z} = \mu(z,w').
\end{equation}
Furthermore, we also have $h_{s,w'}^{w} = 1$.
\end{lem}
\begin{dem}
First, we recall that the polynomials $h_{s,w'}^{z}$ were defined by the isomorphism
\begin{equation}
B_{s} \otimes_{R} B_{w'} \cong  \bigoplus_{z\in \simpleset{w}} h_{s,w'}^{z} B_z.
\end{equation}
According to Soergel's conjecture, this isomorphism implies that
\begin{equation}  \label{equation to categorify one}
\sub{H}_s \sub{H}_{w'} =  \sum_{z\in \simpleset{w}} h_{s,w'}^{z} \sub{H}_{z},
\end{equation}
at the Hecke-algebra level. On the other hand, (\ref{alg kl poly dos})  tells us that
\begin{equation}  \label{equation to categorify two}
\sub{H}_s \sub{H}_{w'} =  \sub{H}_{w} + \sum_{z< w} \mu(z,w') \sub{H}_{z}.
\end{equation}
Therefore, the lemma is proved by comparing the coefficients in (\ref{equation to categorify one}) and (\ref{equation to categorify two}).
\end{dem}

\medskip
According to Lemma \ref{graded module write in simple basis}, in $\mathcal{G}_{\sub{w}}$, we have
\begin{equation}   \label{equation cate one}
[\cell{w}{x}]=  \sum_{z\in \simpleset{w}} \gdn{w}{x}{z} [\simple{w}{z}],
\end{equation}
for all $x \leq w$. Hence, if we apply the homomorphism $\mbox{Res}$ to  (\ref{equation cate one}), and if we use Theorem \ref{teo grothendieck branching rules x greater than xs} and Corollary \ref{corollary branching simples}, then we know that, in $\mathcal{G}_{\sub{w'}}$,

\begin{equation}   \label{equation cate dos}
\sum_{u\in \simpleset{w'}} \left(\sum_{z\in \simpleset{w}}h_{s,u}^{z} \gdn{w}{x}{z} \right)[\simple{w'}{u}]= \left\{
                                  \begin{array}{rl}
                                    v^{-1} [\cell{w'}{x}] + [\cell{w'}{x'}], & \mbox{if } x'<x ; \\
                                     v [\cell{w'}{x}] + [\cell{w'}{x'}], & \mbox{if }  x'>x.
                                  \end{array}
                                \right.
\end{equation}

Furthermore, if we recall that $ \{ [\simple{w'}{u}]\}_{u\in \simpleset{w'}}$ is a $\Z[v,v^{-1}]$-basis of $\mathcal{G}_{\sub{w'}}$, then (\ref{equation cate dos}) provides a family of equations in $\Z[v,v^{-1}]$ (one for each element in $\simpleset{w'}$). In general, there is no simple, direct method for determining $\simpleset{w'}$. However, it is easy to conclude that $w'\in \simpleset{w'}$. The equation obtained by taking the coefficient of $[\simple{w'}{w'}]$ in (\ref{equation cate dos}) is

\begin{equation}   \label{equation cate tres}
\sum_{z\in \simpleset{w}} h_{s,w'}^{z} \gdn{w}{x}{z} = \left\{
                                  \begin{array}{rl}
                                    v^{-1}  \gdn{w'}{x}{w'}  +  \gdn{w'}{x'}{w'}, & \mbox{if } x'<x ; \\
                                     v \gdn{w'}{x}{w'}  +  \gdn{w'}{x'}{w'}, & \mbox{if }  x'>x.
                                  \end{array}
                                \right.
\end{equation}

Consequently, we can apply Proposition \ref{equality graded decompostion numbers with KL} and Lemma \ref{lemma to categorify} to rewrite (\ref{equation cate tres}) as
\begin{equation}   \label{equation cate cuatro}
h_{x,w}+\sum_{ \begin{subarray}{c}
        z\in \simpleset{w} \\ z\neq w
      \end{subarray}}
\mu(z,w') h_{x,z} = \left\{
                                  \begin{array}{rl}
                                    v^{-1}  h_{x,w'}  +  h_{x',w'}, & \mbox{if } x'<x ; \\
                                     v h_{x,w'}  +  h_{x',w'}, & \mbox{if }  x'>x.
                                  \end{array}
                                \right.
\end{equation}

Reordering the terms in this equation gives us

\begin{equation}   \label{equation cate cinco}
\displaystyle h_{x,w} = \left\{ \begin{array}{rl}
                                    v^{-1}  h_{x,w'}  +  h_{x',w'}- \displaystyle\sum_{\substack{z\in \simpleset{w} \\ z\neq w}} \mu(z,w') h_{x,z}, & \mbox{if } x'<x ; \\        &    \\
                                     v h_{x,w'}  +  h_{x',w'}-\displaystyle\sum_{\substack{z\in \simpleset{w} \\ z\neq w}} \mu(z,w') h_{x,z}, & \mbox{if }  x'>x.
                                  \end{array}
                                \right.
\end{equation}

Finally, we note that (\ref{alg kl poly cinco}) and (\ref{equation cate cinco}) are the same equation. Summing up, we have categorified the recursive formula for Kazhdan-Lusztig polynomials given in (\ref{alg kl poly cinco}).





  \bibliographystyle{alpha}


\end{document}